\documentclass[11pt, letterpaper,final]{amsart}
\usepackage{amsfonts,amsmath, amsthm, amssymb, latexsym, epsfig}
\usepackage[english]{babel}
\usepackage[utf8]{inputenc}
\usepackage[all]{xy}
\usepackage{xspace}
\usepackage{comment}
\usepackage{setspace}
\usepackage{enumerate}
\usepackage{stmaryrd}
\usepackage{nicefrac}
\usepackage[mathscr]{eucal}
\usepackage[notcite,notref]{showkeys}

	\evensidemargin \oddsidemargin

\newcommand{\Ad}{\mathrm{Ad}\,}
	
	\theoremstyle{plain}

	\theoremstyle{definition}

	\theoremstyle{plain}

	\newtheorem*{theorem*}{Theorem}
	\theoremstyle{definition}
	\newtheorem*{remark*}{Remark}
	
	\theoremstyle{definition}

	\numberwithin{equation}{section}
	\numberwithin{figure}{section}

\begin{document}
	\title{Classification of quasi-free compact group actions on Cuntz algebras}
 \author{James Gabe}
\address{Department of Mathematics and Computer Science, University of Southern Denmark, 5230 Odense, Denmark}
        \email{gabe@imada.sdu.dk}
        \thanks{This was supported by DFF grants 1054-00094B and 1026-00371B}
        \subjclass[2020]{46L35, 46L55, 46L80, 20C15}

\begin{abstract}
Quasi-free actions of finite groups on Cuntz algebras $\mathcal O_n$ for $n\geq 2$ are classified up to conjugacy by data in the representation ring. Partial results are obtained for quasi-free actions by compact groups. 
\end{abstract}
\maketitle
\thispagestyle{empty}

Let $G$ be a 2\textsuperscript{nd} countable compact group and $\pi \colon G \to \mathcal U(\mathbb C^n)$ a faithful unitary representation for $n\geq 2$. This induces a $G$-action $\alpha^{\pi}$ on the Cuntz algebra $\mathcal O_n$ (the $C^\ast$-algebra generated by isometries $s_1,\dots, s_n$ with $\sum_{j=1}^n s_js_j^\ast = 1$) called a \emph{quasi-free} action, given by $\alpha^{\pi}_g (s_j) = \sum_{k=1}^n \pi(g)_{k,j} s_k$ for $g\in G$.
Such actions and their fixed point algebras (known as Doplicher--Roberts algebras) have been studied for instance in \cite{Evans, DR-Annals, DR-Inventiones, Izumi-Duke04, Izumi-Adv04, Izumi-24}.  Goldstein and Izumi showed \cite{GoldsteinIzumi} that quasi-free outer $G$-actions on $\mathcal O_\infty$  (defined analogously) for $G$ finite are all conjugate, and Izumi recently proved the same for compact groups \cite[Corollary 7.3]{Izumi-24} under the same condition as in the following theorem. 
The following theorem classifies quasi-free actions on $\mathcal O_n$ by using a recent classification theorem \cite{GabeSzabo-Acta}.

In the following, $R(G)$ denotes the (complex) representation ring of $G$. Recall that two elements $x,y\in R(G)$ are \emph{associated} if there is an invertible $u\in R(G)^\times$ such that $ux = y$. Moreover, $\mathcal F(\mathbb C^n)$ denotes the full Fock space $\bigoplus_{k=0}^\infty (\mathbb C^n)^{\otimes k}$, and $\mathcal F(\pi) \colon G \to \mathcal U(\mathcal F(\mathbb C^n))$ denotes the induced representation. Case (1) below is actually a consequence of case (2), see \cite[Chapter XV, Theorem IV]{Burnside}, but this is not needed in the proof.

 \begin{theorem*}
 Let $G$ be a 2\textsuperscript{nd} countable compact group, let $n\geq 2$, and let $\pi_1, \pi_2 \colon G \to \mathcal U(\mathbb C^n)$ be faithful representations. Suppose that (1) $G$ is finite; or (2) every irreducible representation of $G$ is equivalent to a subrepresentation of $\mathcal F(\pi_1)$ and $\mathcal F(\pi_2)$. 
Then $\alpha^{\pi_1}$ and $\alpha^{\pi_2}$ are conjugate (i.e.~there exists $\beta \in \mathrm{Aut}(\mathcal O_n)$ with $\beta  \alpha^{\pi_1} = \alpha^{\pi_2} \beta $) if and only if $1-[\pi_1]$ and $1-[\pi_2]$ are associated in the representation ring $R(G)$. 
 \end{theorem*}
\begin{proof}
As seen for instance in \cite[Section 8]{GoldsteinIzumi}, the $G$-equivariant Cuntz--Toeplitz extension
\[
0 \to (\mathcal K(\mathcal F(\mathbb C^n)),  \Ad \mathcal F(\pi_j)) \to (\mathcal E_n, \Ad \mathcal F(\pi_j)) \to (\mathcal O_n, \alpha^{\pi_j}) \to 0
\]
has a $G$-equivariant completely positive splitting, and hence induces an exact triangle in $KK^G$ (see \cite[Appendix A]{MeyerNest-BaumConnes} for details on $KK^G$ as a triangulated category). By results of Pimsner \cite[Lemma 4.7, Remark 4.10(2)]{Pimsner-CuntzPimsner}, there are canonical $KK^G$-equiva\-lences $(\mathcal K(\mathcal F(\mathbb C^n)), \Ad \mathcal F(\pi_j)) \sim_{KK^G} \mathbb C$ and $(\mathcal E_n, \Ad \mathcal F(\pi_j)) \sim_{KK^G} \mathbb C$ which induce an isomorphism in $KK^G$ of this exact triangle and a triangle of the form
\[
\Sigma  (\mathcal O_n, \alpha^{\pi_j}) \to \mathbb C \xrightarrow{1 - [\pi_j]} \mathbb C \xrightarrow{\iota_{j}}  (\mathcal O_n, \alpha^{\pi_j})
\]
where $\iota_{j} \in KK^G(\mathbb C, (\mathcal O_n, \alpha^{\pi_j}))$ is induced by the unital inclusion $\mathbb C \hookrightarrow \mathcal O_n$, and $1-[\pi_j] \in KK^G(\mathbb C, \mathbb C) \cong R(G)$ (see \cite[Remark 2.15(1)]{Kasparov} for this isomorphism).  
Consider the following diagram in $KK^G$ whose rows are exact triangles
\[
\xymatrix{
\Sigma  (\mathcal O_n, \alpha^{\pi_1})\ar[r] \ar@{.>}[d]^{\Sigma x} & \mathbb C \ar@{.>}[d]^u \ar[r]^{1-[\pi_1]} & \mathbb C \ar@{=}[d] \ar[r]^{\iota_1\qquad} &  (\mathcal O_n, \alpha^{\pi_1}) \ar@{.>}[d]^x  \\
\Sigma  (\mathcal O_n, \alpha^{\pi_2}) \ar[r] & \mathbb C \ar[r]^{1- [\pi_2]} & \mathbb C \ar[r]^{\iota_2 \qquad} &  (\mathcal O_n, \alpha^{\pi_2}).
}
\]
By the axioms of triangulated categories and \cite[Proposition 1.1.20]{Neeman}, there exists a $KK^G$-equivalence $x$ making the right square commute if and only if there exists a $KK^G$-equivalence $u$ making the middle square commute. We note that $\alpha^{\pi_j}$ is isometrically shift-absorbing by \cite[Proposition 3.15]{GabeSzabo-Acta} in case (1) (if $G$ is finite), or by \cite[Example 7.1]{Izumi-24} in case (2). By \cite[Corollary 6.4]{GabeSzabo-Acta}, the existence of an $x$ as above is equivalent to the existence of a conjugacy between $\alpha^{\pi_1}$ and $\alpha^{\pi_2}$, while existence of such a $u \in KK^G(\mathbb C, \mathbb C) \cong R(G)$ simply means that $1-[\pi_1]$ and $1-[\pi_2]$ are associated in $R(G)$.
\end{proof}

A similar result for $C^\ast$-algebras generated by vector bundles was obtained by Dadarlat in \cite[Theorem 1.1]{Dadarlat}. %He refers to the proof technique as a homotopy miracle.

I suspect it is possible to take a conjugacy between $\alpha^{\pi_1}$ and $\alpha^{\pi_2}$ and induce an invertible element in $R(G)$ implementing the associatedness of $1-[\pi_1]$ and $1-[\pi_2]$, but a proof of this eludes me. I strongly doubt that the converse is possible.

In the above theorem, the assumption that $\mathcal F(\pi_j)$ contains every irreducible representation was to ensure that the action $\alpha^{\pi_j}$ is isometrically shift-absorbing. By \cite[Theorem 1.2]{Izumi-24} (one implication is proved in \cite{Mukohara}) this is equivalent to the fixed-point algebra $\mathcal O_n^{\alpha^{\pi_j}}$ (the Doplicher-Roberts algebra) being simple and purely infinite (since $\alpha^{\pi_j}$ is minimal). 
It is suggested in \cite[Example 7.1]{Izumi-24} that the condition in the above theorem might also be necessary for $\alpha^{\pi_j}$ being isometrically shift-absorbing.

I thank Gábor Szabó for useful conversations, Wojciech Szyma\'nski for helpful comments and for pointing out the reference \cite{Dadarlat}, and Masaki Izumi for suggesting generalising the result from finite to compact groups. Finally, I thank the referee for useful comments.

\end{document}